\documentclass[12pt]{article}
\usepackage{amssymb,amsfonts,amsmath,amsthm}

\def\3{\subset }
\def\4{\subseteq }
\def\ov{\overline}

\def\0{\leqno}

\def\barr{\begin{array}}
\def\earr{\end{array}}
\def\dd{\displaystyle}

\def\Z{{\rlap{$\kern2pt{\rm Z}$}{\rm Z}\,}}
\def\bld#1#2{{\buildrel{#1}\over{#2}}}
\def\st#1#2{{\mathrel{\mathop{#2}\limits_{#1}}{}\!}}
\def\stb#1#2#3{{\st{{#1}}{\bld{{#2}}{#3}}{}\!}}
\def\xmare#1#2{\stb{#1}{#2}{\mbox{\Huge$\times$}}}

\title{\bf A note on subgroup coverings\\ of finite groups}
\author{Marius T\u arn\u auceanu}
\date{October 1, 2011}

\begin{document}

\maketitle

\begin{abstract}
In this note we determine the finite groups that can be written as
the union of any three irredundant/distinct proper subgroups. The
finite groups that can uniquely be written as the union of three
proper subgroups are also characterized.
\end{abstract}
\bigskip


\bigskip

\section{Introduction}

A well-known elementary result of group theory states that a group
cannot be written as the union of two proper subgroups. In Scorza
\cite{9} the groups which are the union of three proper subgroups
have been characterized. The analogous problems with three
replaced by four, five and six subgroups were solved by Cohn
\cite{3}, while the case of seven subgroups was studied by
Tomkinson \cite{8}.

Following Cohn's notation, for a group $G$ we shall write
$\sigma(G)=n$ whenever $G$ is the union of $n$ proper subgroups,
but is not the union of any smaller number of proper subgroups. By
using this notation, we first recall the above mentioned results.

\bigskip\noindent{\bf Theorem A.} {\it Let $G$ be a group. Then
\begin{itemize}
\item[\rm a)]
        $\sigma(G) \notin \{1, 2\}${\rm ;}
\item[\rm b)]
        $\sigma(G)=3$ if and only if $G$ has a quotient isomorphic to $C_2 \times C_2${\rm ;}
\item[\rm c)]
        $\sigma(G)=4$ if and only if $\sigma(G) \neq 3$ and $G$ has a quotient isomorphic to $C_3 \times C_3$ or $S_3${\rm ;}
\item[\rm d)]
        $\sigma(G)=5$ if and only if $\sigma(G) \notin \{3, 4\}$ and $G$ has a quotient isomorphic to $A_4${\rm ;}
\item[\rm e)]
        $\sigma(G)=6$ if and only if $\sigma(G) \notin \{3, 4, 5\}$ and $G$ has a quotient isomorphic to $C_5 \times C_5$, $D_5$ or to the group of order {\rm 20} with the following presentation $\langle a,b \mid a^5 = b^4=1, \hspace{1mm}ba=a^2b \rangle${\rm;}
\item[\rm f)]
        $\sigma(G) \neq 7$.
\end{itemize}}

The topic of covering finite groups with proper subgroups has
enjoyed a rapid development in the last few years. Examples of
subjects which have been studied are the following: some special
types of coverings (as minimal coverings, Hall coverings, Sylow
coverings, ... and so on), subgroup co\-ve\-rings of symmetric
groups, alternating groups and linear groups, coverings of groups
with cosets or conjugates of proper subgroups. Many interesting
open problems on this topic can be found in \cite{2}.

The starting point for our discussion is given by the previous
Scorza's result b). A natural question arises from this result:
{\it how many distinct cove-rings with three proper subgroups
possesses a finite group}? In this note we answer partially this
question, by characterizing the finite groups for which any three
irredundant/distinct proper subgroups form a covering (here, three
subgroups are called {\it irredundant} if no one is contained in
the union of the other two) and the finite groups possessing a
unique such covering. Explicit formulas of the number of distinct
coverings with three proper subgroups for several classes of
finite groups are also given.

Most of our notation is standard and will usually not be repeated
here. For basic notions and results on groups we refer the reader
to \cite{1} and \cite{4}.
\bigskip

Clearly, there is no group that can be written as the union of any
three proper subgroups. Therefore some supplementary conditions on
the components of such a covering must be imposed. One of them is
presented in the following result. It is satisfied by two finite
groups: the Klein's group and the quaternion group of order 8.

\bigskip\noindent{\bf Theorem B.} {\it A finite group $G$ is the union of any three irredundant proper subgroups if and only if it is isomorphic to $C_2 \times C_2$ or $Q_8$.}
\bigskip

Another natural condition is obtained by replacing the irredundant
proper subgroups with distinct proper subgroups in Theorem B.
Obviously, if a finite group $G$ is the union of any three
distinct proper subgroups, it is also the union of any three
irredundant proper subgroups. Then $G \cong C_2 \times C_2$ or $G
\cong Q_8$. Since $Q_8$ does not verify our condition (the
collection consisting of any two distinct maximal subgroups and
their intersection is not a covering of $Q_8$), we infer that it
is verified only by $C_2 \times C_2$. Hence the following
corollary holds.

\bigskip\noindent{\bf Corollary C.} {\it A finite group $G$ is the union of any three distinct proper subgroups if and only if it is isomorphic to $C_2 \times C_2$.}
\bigskip

Next, we shall focus on describing the finite groups that can be
covered with three proper subgroups in a unique way.

\bigskip\noindent{\bf Theorem D.} {\it Let $G$ be a finite group. Then the following conditions are equivalent{\rm :}
\begin{itemize}
\item[\rm a)]
        $G$ can uniquely be written as the union of three proper subgroups.
\item[\rm b)]
        $G$ has a unique quotient isomorphic to $C_2 \times C_2$.
\item[\rm c)]
        $G$ has a quotient isomorphic to $C_2 \times C_2$ and no quotient isomorphic to $C_2 \times C_2 \times C_2$.
\end{itemize}}

The condition c) in Theorem D can be rewritten in a more
convenient manner for some particular classes of finite groups.
So, it is easy to see that a finite $p$-group $G$ satisfies c) if
and only if $p=2$ and $(G:\Phi(G))=4$ (that is, $G$ is generated
by exactly two elements). This condition can naturally be extended
to finite nilpotent groups, since such a group is the direct
product of its Sylow subgroups.

\noindent{\bf Corollary E.} {\it Let $G$ be a finite nilpotent
group. Suppose that $G=\xmare{i=1}{k}G_i$, where $G_i$ is the
Sylow $p_i$-subgroup of $G$, $i=\overline{1,k}$, and
$p_1<p_2<...<p_k$. Then $G$ can uniquely be written as the union
of three proper subgroups if and only if $p_1=2$ and $G_1$ is
generated by exactly two elements.}
\bigskip

The above result suggests a way to construct finite groups having
a unique covering with three proper subgroups (for example, the
direct products of types $C_2 \times C_2 \times A$ or $Q_8 \times
A$, where $A$ is an arbitrary finite group of odd order). Also,
the structure of finite hamiltonian groups which satisfy this
property follows immediately from Corollary E.
\bigskip

\noindent{\bf Corollary F.} {\it A finite hamiltonian group can
uniquely be written as the union of three proper subgroups if and
only if it is of type $Q_8 \times A$, where $A$ is a finite
abelian group of odd order.}
\bigskip

Finally, we present a remark that allows us to compute explicitly
the number $c_3(G)$ of distinct coverings with three proper
subgroups of some particular finite groups $G$.
\bigskip

\noindent{\bf Remark.} By the proof of a)$\Longleftrightarrow$ b)
in Theorem D, it will follow that for a finite group there is a
bijection between the set of all coverings with three proper
subgroups and the set of all quotients isomorphic to $C_2 \times
C_2$. In this way, $c_3(G)$ is equal with the number of all
quotients of $G$ which are isomorphic to $C_2 \times C_2$.
\bigskip

Let $G$ be a finite elementary abelian 2-group of order $2^n$.
Then $c_3(G)$ will coincide with the number of all subgroups of
order $2^{n-2}$ in $G$ and can be found by using Corollary 2 of
\cite{6}, \S \hspace{1mm}2.2 (see also \cite{5} and \cite{7}). One
obtains
$$c_3(G)=\frac{2^{2n-1}-3\cdot2^{n-1}+1}{3}\hspace{1mm}.$$

A remarkable class of finite nonabelian groups for which we are
able to compute this number is constituted by the well-known
dihedral groups $$D_{2n}=\langle x,y\mid x^n=y^2=1,
yxy=x^{-1}\rangle, \hspace{1mm}n\geq 2\hspace{1mm}.$$We easily get
that if $n$ is odd $D_{2n}$ possesses no quotient isomorphic to
$C_2 \times C_2$, while if $n$ is even there is only one such
quotient (namely $D_{2n}/\langle x^2\rangle$). Hence
$$c_3(D_{2n})=\left\{\barr{lll}
0,&{\rm if}\hspace{1mm} n \hspace{1mm} {\rm is} \hspace{1mm} {\rm odd}\\
&&\\
1,&{\rm if}\hspace{1mm} n \hspace{1mm} {\rm is} \hspace{1mm} {\rm
even\hspace{1mm}.}\earr\right.$$ Remark that the dihedral groups
$D_{2n}$ with $n\equiv 0 \hspace{1mm}({\rm mod} \hspace{1mm}2)$
are also examples of finite groups having a unique covering with
three proper subgroups.
\bigskip

We finish our note by indicating several open problems concerning
to subgroup coverings of finite groups, which are derived from the
above study.

\bigskip\noindent{\bf Problem 1.} {\rm Give explicit formulas of the number of distinct
coverings with three proper subgroups for other classes of finite
groups.}

\bigskip\noindent{\bf Problem 2.} {\rm Find characterizations of finite groups that
can be written as the union of any three irredundant/distinct
proper normal subgroups. Study the uniqueness of such a covering,
too.}

\bigskip\noindent{\bf Problem 3.} {\rm Study all above problems for coverings of
finite groups with more than three proper (normal) subgroups.}
\bigskip

\section{Proofs of our results}

\bigskip\noindent{\bf Proof of Theorem B.} Obviously, $C_2 \times
C_2$ and $Q_8$ have each of them only three irredundant proper
subgroups and they are the union of these subgroups.

Conversely, let us suppose that $G$ can be written as the union of
any three irredundant proper subgroups and it is not isomorphic to
$C_2 \times C_2$. Then every subgroup of order $p^2$ of $G$ is
cyclic, where $p$ is any prime. In other words, $G$ satisfies the
$p^2 \hspace{1mm} conditions$. It is well-known that these are
equivalent to the $Sylow \hspace{1mm} conditions$, i.e. for any
odd prime $p$ the Sylow $p$-subgroups of $G$ are cyclic and the
Sylow 2-subgroups of $G$ are either cyclic or isomorphic to a
generalized quaternion group $$Q_{2^n} = \langle a,b \mid
a^{2^{n-1}} = b^4 = 1, bab^{-1}=a^{2^{n-1}-1}\rangle,
\hspace{1mm}n \geq 3\hspace{1mm}.$$

Next, let $p$ be an odd prime and $n_p$ be the number of Sylow
$p$-subgroups of $G$. Assume that $n_p \neq 1$. Then there is a
Sylow $p$-subgroup $H$ of $G$ which is not normal. Since $n_p
\equiv 1$ (mod $p$), we infer $n_p \geq p+1 \geq 4$ and so $H$
possesses at least four distinct conjugates, say $H_0=H$, $H_1$,
$H_2$ and $H_3$. By our hypothesis one obtains $$G =
\dd\bigcup_{i=1}^3 H_i,$$which leads to $$H = \dd\bigcup_{i=1}^3
\hspace{1mm}(H \cap H_i).$$Because $H$ is cyclic, we must have $H
= H \cap H_i$ for some $i \in \{1, 2, 3\}$. This shows that $H =
H_i$, a contradiction. Thus every Sylow $p$-subgroup of $G$ is
normal, for any odd prime $p$.

The above conclusion implies that $G=G'G''$, where $G'$ is a
2-group which is cyclic or isomorphic to $Q_{2^n}$ and $G''$ is a
cyclic group of odd order. Moreover, we have $G'' \lhd G$. Now, we
easily obtain that $G'$ cannot be cyclic, since $G$ has a quotient
isomorphic to $C_2 \times C_2$. Then $G' \cong Q_{2^n}$ and
therefore $G'$ can be written as the union of three irredundant
proper subgroups (for example, the union of its maximal
subgroups). This leads to $G = G'$, in view of our hypothesis.

A simple calculation shows that $Q_{2^n}$ possesses $2^{n-2}+1$
distinct subgroups of order 4, which are all cyclic. For $n \geq
4$ we have $2^{n-2}+1 \geq 5$, and it is clear that the union of
three such subgroups is not equal with $Q_{2^n}$. Hence $n=3$ and
$G \cong Q_8$. \hfill\rule{1,5mm}{1,5mm}
\bigskip

\bigskip\noindent{\bf Proof of Theorem D.} We shall prove first a) $\Longleftrightarrow$ b) and then b) $\Longleftrightarrow$ c).

a) $\,\Rightarrow$ b). Suppose that there are two distinct normal
subgroups $H$ and $K$ of $G$ such that $G / H \cong G / K \cong
C_2 \times C_2$. Then $G$ can be written as the following unions
of irredundant proper subgroups $$G = \dd\bigcup_{i=1}^3 H_i =
\dd\bigcup_{i=1}^3 K_i,$$ where $(H_i)_{i=\overline{1,3}}$ and
$(K_i)_{i=\overline{1,3}}$ are the maximal subgroups of $G$ which
contain $H$ and $K$, respectively. This contradicts a), since by
$H \neq K$ we clearly have $\{H_1, H_2, H_3\} \neq \{K_1, K_2,
K_3\}$.

b) $\,\Rightarrow$ a). Let $(H_i)_{i=\overline{1,3}}$ be a
collection of three irredundant proper subgroups of $G$ such that
$$G = \dd\bigcup_{i=1}^3 H_i$$and denote by $H$ the intersection of $H_i, i=\ov{1,3}.$

First of all, we show that $H$ is normal in $G$. Let $h \in H$, $x
\in G$ be two arbitrary elements and assume that $xhx^{-1}
\hspace{-1mm}\notin \hspace{-1mm}H_1$. Then $xhx^{-1} \in
H_2\setminus (H_1 \cup H_3)$ or $xhx^{-1} \in H_3\setminus (H_1
\cup H_2)$. Assume now that $xhx^{-1} \in H_2\setminus (H_1 \cup
H_3)$ (the other situation being analogous). One obtains $x \notin
H_1 \cup H_3$ and so $x$ is contained in $H_2\setminus (H_1 \cup
H_3)$, too. Take an element $y \in H_1\setminus (H_2 \cup H_3)$.
Then both $xhy^{-1}$ and $yx^{-1}$ belong to $H_3\setminus (H_1
\cup H_2)$. It follows that $xhx^{-1} = xhy^{-1}yx^{-1}\in H_3$, a
contradiction. In this way, $xhx^{-1}$ is contained in each $H_i$,
$i=1,2,3$, and therefore in $H$. Thus $H$ is normal in $G$.

A well-known exercise of group theory shows that $x^2 \in H$, for
all $x \in G$. This implies that the quotient $G / H$ is in fact
an elementary abelian 2-group. By our hypothesis, we infer that $G
/ H \cong C_2 \times C_2$. Thus we have proved that any
decomposition of $G$ as a union of three irredundant proper
subgroups induces a quotient of $G$ isomorphic to $C_2 \times
C_2$. Because $G$ possesses only one such quotient, it can
uniquely be covered with three irredundant proper subgroups, as
claimed.

b) $\,\Rightarrow$ c). It is obvious, since $C_2 \times C_2 \times
C_2$ has more quotients isomorphic to $C_2 \times C_2$ (in fact,
exactly seven).

c) $\,\Rightarrow$ b). Suppose that $G$ has two distinct quotients
isomorphic to $C_2 \times C_2$, say $G / H$ and $G / K$. We obtain
$(G:H \cap K) \in \{2^3, 2^4\}$. On the other hand, we know that
$x^2 \in H \cap K$, for all $x \in G$, which implies that $G / H
\cap K$ is elementary abelian. So, it is isomorphic to $C_2 \times
C_2 \times C_2$ or $C_2 \times C_2 \times C_2 \times C_2$, in both
cases contradicting the hypothesis. Hence b) holds.
\hfill\rule{1,5mm}{1,5mm}

\vspace*{5ex}\small

\hfill
\begin{minipage}[t]{5cm}
Marius T\u arn\u auceanu \\
Faculty of  Mathematics \\
``Al.I. Cuza'' University \\
Ia\c si, Romania \\
e-mail: {\tt tarnauc@uaic.ro}
\end{minipage}

\end{document}